\documentclass[11pt]{amsart}
\usepackage{mathrsfs}
\usepackage{amssymb}

\usepackage{titletoc}
\usepackage{amscd}
\usepackage{amsmath, amssymb}
\usepackage{amsfonts}
\usepackage[colorlinks,linkcolor=blue,citecolor=blue, pdfstartview=FitH]{hyperref}

\numberwithin{equation}{section}

\theoremstyle{plain}
\newtheorem{thm}{Theorem}[section]
\newtheorem{theorem}[thm]{Theorem}
\newtheorem{lemma}[thm]{Lemma}
\newtheorem{corollary}[thm]{Corollary}
\newtheorem{proposition}[thm]{Proposition}

\theoremstyle{definition}

\newtheorem{remark}[thm]{Remark}

\newtheorem{definition}[thm]{Definition}

\newtheorem{example}[thm]{Example}

\newtheorem{defn-thm}[thm]{Definition-Theorem}

\renewcommand{\S}{{\mathbb S}}

\newcommand{\bp}{\bar{\partial}}

\newcommand{\btheorem}{\begin{theorem}}
\newcommand{\etheorem}{\end{theorem}}
\newcommand{\bproposition}{\begin{proposition}}
\newcommand{\eproposition}{\end{proposition}}
\newcommand{\bdefinition}{\begin{definition}}
\newcommand{\edefinition}{\end{definition}}
\newcommand{\bcorollary}{\begin{corollary}}
\newcommand{\ecorollary}{\end{corollary}}
\newcommand{\bproof}{\begin{proof}}
\newcommand{\eproof}{\end{proof}}
\newcommand{\bremark}{\begin{remark}}
\newcommand{\eremark}{\end{remark}}
\newcommand{\eexample}{\end{example}}
\newcommand{\bexample}{\begin{example}}

\newcommand{\elemma}{\end{lemma}}
\newcommand{\blemma}{\begin{lemma}}

\newcommand{\sq}{\sqrt{-1}}

\newcommand{\p}{\partial}

\renewcommand{\bar}{\overline}

\renewcommand{\phi}{\varphi}

\newcommand{\ee}{\end{eqnarray*}}
\newcommand{\be}{\begin{eqnarray*}}

\newcommand{\beq}{\begin{equation}}
\newcommand{\eeq}{\end{equation}}

\newcommand{\bd}{\begin{enumerate}}
\newcommand{\ed}{\end{enumerate}}

\newcommand{\Ric}{\mathrm {Ric}}

\newcommand{\ov}[1]{\overline{#1}}
\newcommand{\mn}{\sqrt{-1}}
\newcommand{\de}{\partial}

\usepackage{fancyhdr}
\renewcommand{\leftmark}{{\scriptsize The Chern-Ricci flow and holomorphic bisectional curvature}}

\begin{document}
\title{The Chern-Ricci flow and holomorphic bisectional curvature}
\makeatletter
\let\uppercasenonmath\@gobble% disables title uppercase
\let\MakeUppercase\relax% disables author uppercase
\let\scshape\relax% disables section smallcaps
\makeatother
\author{ Xiaokui Yang}
\date{}
\address{\textsf{Current address of Xiaokui Yang: Morningside Center of Mathematics, Institute of
        Mathematics, Hua Loo-Keng Key Laboratory of Mathematics,
        Academy of Mathematics and Systems Science,
        Chinese Academy of Sciences, Beijing, 100190, China.}}
\email{\href{mailto:xkyang@amss.ac.cn}{\texttt{xkyang@amss.ac.cn}}}

\begin{abstract} In this note, we show that on Hopf manifold $\S^{2n-1}\times \S^1$,  the non-negativity of the holomorphic bisectional curvature is not preserved along the Chern-Ricci flow.

\end{abstract}
\maketitle

\section{Introduction}

The Chern-Ricci flow is an evolution equation for Hermitian metrics
on complex manifolds, generalizing the K\"ahler-Ricci flow.  Given
an initial  Hermitian metric $\omega_0 = \sqrt{-1}
(g_0)_{i\ov{j}}dz^i \wedge d\ov{z}^j$, the Chern-Ricci flow is
defined as
\begin{equation} \label{crf0}
\frac{\p \omega}{\de t}=-\Ric(\omega), \qquad \omega|_{t=0}
=\omega_0,
\end{equation}
where $\Ric(\omega) := - \sq\p\bp \log \det g$ is the Chern-Ricci
form of $\omega$.  In the case when $\omega_0$ is K\"ahler, namely
$d\omega_0=0$, (\ref{crf0}) coincides with the K\"ahler-Ricci flow.
 The Chern-Ricci flow was first introduced by Gill \cite{G} in the setting of manifolds with vanishing first Bott-Chern classes, and
 many
 fundamental properties are established by Tosatti and Weinkove
 \cite{TW} on more general manifolds.   A variety of further results on Chern-Ricci flow are studied in
\cite{TW,TW2, TWY, Gill, GS,  FTWZ,Z} and some of them are analogues
to classical results for the K\"ahler-Ricci flow (e.g.\cite{Ha, Ca,
ST, ST2, TZha, SW1,SW}).

  It is proved by Mok \cite{Mok} (see \cite{Bando} for K\"ahler threefolds and also \cite{Gu}) that the non-negativity of the
  holomorphic bisectional curvature is preserved along the K\"ahler-Ricci
  flow. However, we show that on Hermitian manifolds, the
  non-negativity of the
  holomorphic bisectional curvature is not necessarily preserved
  under the Chern-Ricci flow.

\btheorem\label{main}  Let $X=\S^{2n-1}\times \S^1$ be a diagonal
Hopf manifold. Fix $T_0\geq 0$ and let
$$ \omega_0=\frac{1}{|z|^4}\sum\left((1+T_0)\delta_{ij}|z|^2-T_0 \bar
z^i z^j\right) \sq dz^i\wedge d\bar z^j.$$ Then the Chern-Ricci flow
(\ref{crf0})
 has %explicit solution \beq
%\omega(t)= \omega_0-t\Ric(\omega_0) \qtq{for} 0\leq
%t<\frac{T_0+1}{n},\eeq where $\frac{T_0+1}{n}$ is the
maximal existence time $ T_{\max}=\frac{T_0+1}{n}$. \bd\item When $
t\in\displaystyle\left[0,\frac{T_0}{n}\right]$, $\omega(t)$ has
non-negative holomorphic bisectional curvature;
 \item However, when $
t\in\displaystyle\left(\frac{2T_0+1}{2n},\frac{T_0+1}{n}\right)$,
the holomorphic bisectional curvature of $\omega(t)$ is no longer
non-negative. \ed \etheorem \bremark It is worth to point out that
the same proof as in the K\"ahler case (following Mok) fails for the
Chern-Ricci flow  since the evolution of the Riemann curvature
tensor under the Chern-Ricci flow involves also some terms with the
torsion (and its covariant derivatives), which are not there in the
K\"ahler-Ricci flow, where the evolution of the curvature involves
only the curvature tensor itself.

\eremark

\bremark It is also interesting to investigate sufficient conditions
on Hermitian manifolds such that the non-negativity of the
holomorphic bisectional curvature is preserved under the Chern-Ricci
flow.

\eremark

\section{The proof of Theorem \ref{main}}

For $\alpha = (\alpha_1, \ldots, \alpha_n) \in \mathbb{C}^n
\setminus \{ 0 \}$ with $|\alpha_1|=\dots=|\alpha_n| \neq 1$, let
$M$ be the Hopf manifold
 $M=(\mathbb{C}^n\setminus \{0\})/\sim$, where
$$(z^1, \dots,z^n) \sim \left(\alpha_1 z^1, \dots, \alpha_n
z^n\right).$$ It is easy to see that  $M$ is diffeomorphic to
$\S^{2n-1}\times \S^1$.  Fix $T_0>0$ and  consider the Hermitian
metric
$$ \omega_0=\frac{1}{|z|^4}\left((1+T_0)\delta_{ij}|z|^2-T_0 \bar
z^i z^j\right) \sq dz^i\wedge d\bar z^j.$$ where $|z|^2=\sum_{j=1}^n
|z^j|^2$. It is proved in \cite{TW} that \beq \omega(t) = \omega_0 -
t \Ric(\omega_0)\eeq gives an explicit solution of the Chern-Ricci
flow on $M$ with initial metric $\omega_0$. Indeed, by elementary
linear algebra, we see $\det(\omega_0)=(1+T_0)^{n-1}|z|^{-2n}$ and
so
$$\Ric(\omega_0) = n \sq \p\bp \log |z|^2 = \frac{n}{|z|^2} \left(\delta_{ij} - \frac{\ov{z}^i z^j}{|z|^2}\right)\mn dz^i \wedge d\ov{z}^j \ge 0.$$
For $t<\frac{T_0+1}{n}$, we have the Hermitian metrics \beq
\omega(t)=\omega_0-t\Ric(\omega_0)=\frac{1}{|z|^2}\left((1+T_0-nt)\delta_{ij}-(T_0-nt)\frac{\ov{z}^i
z^j}{|z|^2}\right) \mn dz^i \wedge d\ov{z}^j.\label{time}\eeq Hence
$$\det(\omega(t))=\frac{(1+T_0-nt)^{n-1}}{|z|^{2n}},$$
from which it follows that
$\Ric(\omega(t))=\Ric(\omega_0)=n\sq\p\bp\log|z|^2.$  It also
implies that $\omega(t)$ solves the Chern-Ricci flow on the maximal
existence interval $\displaystyle\left[0,\frac{T_0+1}{n}\right).$

  Next, we compute the curvature tensor of the involving metric
  (\ref{time}). For simplicity, we define a rescaled metric
$\omega_\lambda=\sq h_{i\bar j}dz^i\wedge d\bar z^j$  on $M$ with
\beq h_{i\bar j}=\frac{1}{|z|^4}\left(\delta_{ij}|z|^2-\lambda \bar
z^i z^j\right),\ \ \ \ \ \ \ \lambda<1. \eeq Note that when
$$\lambda=\frac{T_0-nt}{1+T_0-nt},$$ we have \beq
\omega_\lambda=\frac{\omega(t)}{1+T_0-nt}.\label{relation}\eeq

\blemma\label{formula} Let $R_{ k\bar j i\bar q}$ be the curvature
components of $\omega_\lambda$, then  \be R_{k\bar j i \bar q}
&=&\frac{\delta_{iq}(\delta_{jk}|z|^2-\bar z^k
z^j)}{|z|^6}+\frac{\lambda \left(\delta_{ij}|z|^2-\bar z^i
z^j\right)\left(\delta_{kq}|z|^2-\bar z^k
z^q\right)}{|z|^8}\\&&+\frac{(\lambda^2-2\lambda)\bar
z^iz^q(\delta_{kj}|z|^2-\bar z^k z^j)}{|z|^8}. \ee \elemma

 \bproof By using elementary linear algebra, one has
$\det(h_{i\bar j})=(1-\lambda)|z|^{-2n}$ and so \beq
\Ric(\omega_\lambda)=n\sq \p\bp\log |z|^2\geq 0.\eeq
 On the other hand, one can verify that the matrix $(h_{i\bar j})$ has (transpose) inverse matrix \beq h^{i\bar
j}=|z|^2\left(\delta_{ij}+\frac{\lambda z^i\bar
z^j}{(1-\lambda)|z|^2}\right).\eeq  By straightforward computation,
\beq \frac{\p h_{i\bar j}}{\p z^k}=-\frac{\delta_{ij}\bar
z^k}{|z|^4}-\frac{\lambda\delta_{jk}\bar z^i}{|z|^4}+\frac{2\lambda
\bar z^i\bar z^k z^j}{|z|^6}=\frac{2\lambda \bar z^i\bar z^k
z^j}{|z|^6}-\frac{\lambda \delta_{jk}\bar z^i+\delta_{ij}\bar
z^k}{|z|^4}\eeq
 and so
 \be \Gamma_{ki}^p&=&h^{p\bar j}\frac{\p h_{i\bar
j}}{\p z^k}=|z|^2\left(\delta_{pj}+\frac{\lambda z^p\bar
z^j}{(1-\lambda)|z|^2}\right)\left(\frac{2\lambda \bar z^i\bar z^k
z^j}{|z|^6}-\frac{\lambda \delta_{jk}\bar z^i+\delta_{ij}\bar
z^k}{|z|^4}\right)\\
&=&\frac{2\lambda\bar z^i\bar z^kz^p}{|z|^4}-\frac{\lambda
\delta_{pk}\bar z^i+\delta_{ip}\bar z^k}{|z|^2}
+\frac{2\lambda^2\bar z^i\bar z^k
z^p}{(1-\lambda)|z|^4}-\frac{\lambda^2\bar z^i\bar z^k z^p+\lambda
\bar z^i\bar z^k z^p}{(1-\lambda)|z|^4}\\
&=&\frac{\lambda\bar z^i\bar z^kz^p}{|z|^4}-\frac{\lambda
\delta_{pk}\bar z^i+\delta_{ip}\bar z^k}{|z|^2}.\ee

\noindent The Chern curvature tensor of $\omega_\lambda$ is \be
&&R_{k\bar j i}^p=-\frac{\p\Gamma_{ki}^p}{\p\bar
z^j}\\&=&-\frac{\lambda \delta_{ij}\bar z^k z^p+\lambda
\delta_{kj}\bar z^i z^p}{|z|^4}+\frac{2\lambda \bar z^i\bar z^k z^p
z^j}{|z|^6}+\frac{\lambda
\delta_{pk}\delta_{ij}+\delta_{ip}\delta_{kj}}{|z|^2}-\frac{\lambda
\delta_{pk}\bar z^i z^j+\delta_{ip}\bar z^k z^j}{|z|^4}\\
&=&\frac{\lambda
\delta_{pk}\delta_{ij}+\delta_{ip}\delta_{kj}}{|z|^2}+\frac{2\lambda
\bar z^i\bar z^k z^p z^j}{|z|^6}-\frac{\lambda\left(\delta_{ij}\bar
z^k z^p+ \delta_{kj}\bar z^i z^p+\delta_{pk}\bar z^i
z^j\right)+\delta_{ip}\bar z^k z^j}{|z|^4}. \ee

\noindent Hence \small{\be &&R_{k\bar j i \bar q}=h_{p\bar q}
R_{k\bar j i}^p\\&=&\frac{\delta_{pq}|z|^2}{|z|^4}\left[
\frac{\lambda
\delta_{pk}\delta_{ij}+\delta_{ip}\delta_{kj}}{|z|^2}+\frac{2\lambda
\bar z^i\bar z^k z^p z^j}{|z|^6}-\frac{\lambda\left(\delta_{ij}\bar
z^k z^p+ \delta_{kj}\bar z^i z^p+\delta_{pk}\bar z^i
z^j\right)+\delta_{ip}\bar z^k z^j}{|z|^4}\right]\\&&-\frac{\lambda
\bar z^p z^q}{|z|^4}\left[ \frac{\lambda
\delta_{pk}\delta_{ij}+\delta_{ip}\delta_{kj}}{|z|^2}+\frac{2\lambda
\bar z^i\bar z^k z^p z^j}{|z|^6}-\frac{\lambda\left(\delta_{ij}\bar
z^k z^p+ \delta_{kj}\bar z^i z^p+\delta_{pk}\bar z^i
z^j\right)+\delta_{ip}\bar z^k z^j}{|z|^4}\right]\\
&=&\frac{\lambda
\delta_{qk}\delta_{ij}+\delta_{iq}\delta_{jk}}{|z|^4}+\frac{2\lambda
\bar z^i\bar z^k z^j z^q}{|z|^8}-\frac{\lambda\left(\delta_{ij}\bar
z^k z^q+\delta_{kj}\bar z^iz^q+\delta_{kq}\bar z^i
z^j\right)+\delta_{iq}\bar z^k z^j}{|z|^6}\\
&&-\frac{\lambda^2\delta_{ij}\bar z^k z^q+\lambda\delta_{kj}\bar
z^iz^q}{|z|^6}-\frac{2\lambda^2\bar z^i\bar z^k z^j
z^q}{|z|^8}\\&&+\frac{\lambda^2\left(\delta_{ij}\bar z^k
z^q|z|^2+\delta_{kj}\bar z^i z^q|z|^2+\bar z^i\bar z^k z^j
z^q\right)+\lambda \bar z^i\bar z^k z^j z^q}{|z|^8}\\
&=&\frac{\lambda
\delta_{qk}\delta_{ij}+\delta_{iq}\delta_{jk}}{|z|^4}+\frac{(3\lambda-\lambda^2)\bar
z^i\bar z^k z^j z^q}{|z|^8}-\frac{\lambda \delta_{qk}\bar z^i
z^j}{|z|^6}-\frac{\lambda\delta_{ij}\bar z^k
z^q}{|z|^6}\\&&+\frac{(\lambda^2-2\lambda)\delta_{kj}\bar z^i
z^q}{|z|^6}+\frac{\delta_{iq}\bar z^k z^j}{|z|^6}\\
&=&\frac{\delta_{iq}(\delta_{jk}|z|^2-\bar z^k
z^j)}{|z|^6}+\frac{\lambda \delta_{ij}\left(\delta_{kq}|z|^2-\bar
z^k z^q\right)}{|z|^6}+\frac{(\lambda^2-2\lambda)\bar
z^iz^q(\delta_{kj}|z|^2-\bar z^k z^j)}{|z|^8}\\&&+\frac{\lambda\bar
z^i
z^j\left(\bar z^k z^q-\delta_{kq}|z|^2\right)}{|z|^8}\\
&=&\frac{\delta_{iq}(\delta_{jk}|z|^2-\bar z^k
z^j)}{|z|^6}+\frac{\lambda \left(\delta_{ij}|z|^2-\bar z^i
z^j\right)\left(\delta_{kq}|z|^2-\bar z^k
z^q\right)}{|z|^8}+\frac{(\lambda^2-2\lambda)\bar
z^iz^q(\delta_{kj}|z|^2-\bar z^k z^j)}{|z|^8}. \ee}\eproof

\blemma\label{hbsc}  For any $\lambda \in [0,1)$, $\omega_\lambda$
has non-negative holomorphic bisectional curvature. \elemma \bproof
For any $\xi=(\xi^1,\cdots, \xi^n)$ and $\eta=(\eta^1,\cdots,
\eta^n)$, by Lemma \ref{formula} we have \be R_{k\bar j i\bar q}
\xi^k\bar \xi^j \eta^i\bar
\eta^q&=&\frac{|\eta|^2(|z|^2|\xi|^2-|\bar z\cdot
\xi|^2)}{|z|^6}+\frac{\lambda \left|\left(\delta_{ij}|z|^2-\bar z^i
z^j\right)\eta^i\bar\xi^j\right|^2}{|z|^8}\\&&+\frac{(\lambda^2-2\lambda)|\bar
z \cdot \eta|^2(|z|^2|\xi|^2-|z\cdot\bar \xi|^2)}{|z|^8}.\ee Since
$|z|^2|\eta|^2\geq |\bar z\cdot \eta|^2$, we obtain \be R_{k\bar j
i\bar q} \xi^k\bar \xi^j \eta^i\bar \eta^q\geq \frac{\lambda
\left|\left(\delta_{ij}|z|^2-\bar z^i
z^j\right)\eta^i\bar\xi^j\right|^2}{|z|^8}+\frac{(\lambda^2-2\lambda+1)|\bar
z \cdot \eta|^2(|z|^2|\xi|^2-|z\cdot\bar \xi|^2)}{|z|^8}.\ee The
right hand side is non-negative when $\lambda\geq 0$.\eproof

\bcorollary The initial metric  $\omega_0$ has non-negative
holomorphic bisectional curvature. \ecorollary

\bproof When $t=0$,  or equivalently $\lambda=\frac{T_0}{1+T_0} $,
we know $\omega_\lambda=\frac{\omega_0}{1+T_0}$. Since
$\lambda=\frac{T_0}{1+T_0}\in [0,1)$, by Lemma \ref{hbsc},
$\omega_0$ has non-negative holomorphic bisectional curvature.
\eproof

\blemma\label{no} When $\lambda<-1$, the holomorphic sectional
curvature of the metric $\omega_\lambda$ is no longer non-negative.
In particular, the holomorphic bisectional curvature of the metric
$\omega_\lambda$ is no longer non-negative. \elemma \bproof For any
$\xi=(\xi^1,\cdots, \xi^n)$, we have \be R_{k\bar j i\bar q}
\xi^k\bar\xi^j\xi^i\bar \xi^q&=&\frac{|\xi|^2(|z|^2|\xi|^2-|\bar
z\cdot \xi|^2)}{|z|^6}+\frac{\lambda (|z|^2|\xi|^2-|\bar z\cdot
\xi|^2)^2}{|z|^8}\\&&+\frac{(\lambda^2-2\lambda)|\bar z\cdot
\xi|^2(|z|^2|\xi|^2-|\bar z\cdot \xi|^2)}{|z|^8}\\&=&
\frac{(3\lambda-\lambda^2) |\bar z\cdot
\xi|^4+(\lambda+1)(|z|^2|\xi|^2)^2+(\lambda^2-4\lambda-1) |\bar
z\cdot \xi|^2|z|^2\cdot|\xi|^2}{|z|^8}.\ee

\noindent
 Let $a=|\bar z\cdot
\xi|^2$ and $b=|z|^2|\xi|^2$, then \be R_{k\bar j i\bar q}
\xi^k\bar\xi^j\xi^i\bar
\xi^q&=&\frac{(3\lambda-\lambda^2)a^2+(\lambda^2-4\lambda-1)ab+(\lambda+1)b^2}{|z|^8}\\
&=&\frac{(b-a)a(\lambda-1)^2+(b-a)^2(\lambda+1)}{|z|^8}.\ee It is
easy to see that, $b\geq a\geq 0$ and so for any $-1\leq \lambda<1$
$$R_{k\bar j i\bar q}
\xi^k\bar\xi^j\xi^i\bar \xi^q\geq 0.$$
 However, when
$\lambda<-1$, $R_{k\bar j i\bar q} \xi^k\bar\xi^j\xi^i\bar \xi^q$ is
no longer nonnegative. Indeed, for any given $z=(z^1,\cdots, z^n)$,
we choose a nonzero vector $\xi=(\xi^1,\cdots, \xi^n)$ such that
$\bar z \cdot \xi=0$, i.e. $\sum \bar z^i\cdot \xi^i=0$. In this
case, we have $a=|\bar z\cdot \xi|=0$, but $b=|z|^2|\xi|^2>0.$
Moreover,
$$R_{k\bar j i\bar q} \xi^k\bar\xi^j\xi^i\bar
\xi^q=\frac{b^2(\lambda+1)}{|z|^8}<0$$ since $\lambda<-1$. \eproof

\noindent\emph{The proof of Theorem \ref{main}.} By
(\ref{relation}), we see when $\lambda=\frac{T_0-nt}{1+T_0-nt}$,
$\omega_\lambda=\frac{\omega(t)}{1+T_0-nt}$. Hence,

\bd \item by Lemma \ref{hbsc}, when $\lambda\in [0,1)$ or
equivalently, $0\leq t\leq \frac{T_0}{n}$, $\omega(t)$ has
non-negative holomorphic bisectional curvature;

\item by Lemma \ref{no}, when $\lambda<-1$, or equivalently,
$\frac{2T_0+1}{2n}<t<\frac{T_0+1}{n}$, the holomorphic bisectional
curvature of  $\omega(t)$ is no longer non-negative.\qed \ed

\vskip 2cm

\textbf{Acknowledgements.} This work was partially supported by
China's Recruitment
 Program of Global Experts,
 as well as National Center for Mathematics and Interdisciplinary Sciences,
 Chinese Academy of Sciences.  The author wishes to thank Professor Ben Weinkove
and Professor Valentino Tosatti for their invaluable support and
guidance.  The author is grateful to the referees for  very careful
reading and many helpful suggestions.

\end{document}